\documentclass[12pt]{elsarticle}

\usepackage{amsmath,graphicx}
\usepackage{amssymb}
\usepackage{mathtools}
\usepackage{mathrsfs}
\usepackage[margin=2.5cm]{geometry}

\DeclarePairedDelimiter\ceil{\lceil}{\rceil}

\newcommand{\Le}{{\mathscr L}}

\begin{document}

\begin{frontmatter}


\title{Geodetic convexity and Kneser graphs}

\author[infes-uff]{Marcos Bedo}
\ead{marcosbedo@id.uff.br}

\author[infes-uff]{Jo\~ao V. S. Leite}
\ead{joaovitorleite@id.uff.br}

\author[infes-uff]{Rodolfo A. Oliveira}
\ead{rodolfooliveira@id.uff.br}

\author[ic-uff]{F\'abio Protti\corref{cor1}}
\ead{fabio@ic.uff.br}

\cortext[cor1]{Corresponding author}

\address[infes-uff]{Instituto do Noroeste Fluminense, Universidade Federal Fluminense, Brazil}

\address[ic-uff]{Instituto de Computa\c c\~ao, Universidade Federal Fluminense, Brazil}

\begin{abstract}
The {\em Kneser graph} $K(2n+k,n)$, for positive integers $n$ and $k$, is the graph $G=(V,E)$ such that $V=\{S\subseteq\{1,\ldots,2n+k\} : |S|=n\}$ and there is an edge $uv\in E$ whenever $u\cap v=\emptyset$. Kneser graphs have a nice combinatorial structure, and many parameters have been determined for them, such as the diameter, the chromatic number, the independence number, and, recently, the hull number (in the context of $P_3$-convexity). However, the determination of geodetic convexity parameters in Kneser graphs still remained open. In this work, we investigate both the geodetic number and the geodetic hull number of Kneser graphs. We give upper bounds and determine the exact value of these parameters for Kneser graphs of diameter two (which form a nontrivial subfamily). We prove that the geodetic hull number of a Kneser graph of diameter two is two, except for $K(5,2)$, $K(6,2)$, and $K(8,2)$, which have geodetic hull number three. We also contribute to the knowledge on Kneser graphs by presenting a characterization of endpoints of diametral paths in $K(2n+k,n)$, used as a tool for obtaining some of the main results in this work.
\end{abstract}


\end{frontmatter}

\journal{Applied Mathematics and Computation}

\newtheorem{theorem}{Theorem}
\newtheorem{lemma}[theorem]{Lemma}
\newtheorem{corollary}[theorem]{Corollary}
\newtheorem{proposition}[theorem]{Proposition}
\newtheorem{fact}[theorem]{Fact}
\newtheorem{example}[theorem]{Example}

\section{Introduction}
\label{sec:intro}

Let $n,k$ be positive integers. The {\em Kneser graph} $K(2n+k,n)$ is the graph $G=(V,E)$ such that $V=\{S\subseteq\{1,\ldots,2n+k\}: |S|=n\}$ and there is an edge $uv\in E$ whenever $u\cap v=\emptyset$. Kneser graphs have a rich combinatorial structure~\cite{GODSIL,LOVASZ78}, and there are many studies on this class involving colorings, independent sets, and products of graphs (see~\cite{BALOGH2019,BRESAR2019,JIN2020}). In addition, the $P_3$-hull number of a Kneser graph has been investigated in~\cite{GRIPPO21}, where the authors determine the exact value of the $P_3$-hull number of $K(2n+k,n)$ for $k>1$, and provide lower and upper bounds for $k=1$. Similarly, the authors in~\cite{LIAO2023} have recently studied the $q$-analogues of Kneser graphs under the same convexity and parameter. To the best of the authors' knowledge, however, no studies on geodetic convexity parameters are known for Kneser graphs. This work investigates two of the most addressed  parameters in graph convexity, the {\em geodetic number} and the {\em geodetic hull number}, in the context of Kneser graphs. In particular, we determine the exact value of these parameters for Kneser graphs of diameter two. This turns out to be a relevant question because, for a fixed $n\geq 2$, almost all graphs in the family ${\mathscr K}_n=\{K(2n+k,n) : k\geq 1\}$ have diameter two (those for which $k\geq n-1$). We prove that the geodetic hull number of a Kneser graph of diameter two is two, except for $K(5,2)$, $K(6,2)$, and $K(8,2)$, which have geodetic hull number three.

The geodetic and geodetic hull numbers have been studied for several graph classes, e.g.~\cite{BRESAR20084044,BRESAR20111693,CACERES20103020,GMnumber}. For general graphs, determining the geodetic or the geodetic hull number is NP-hard~\cite{DouradoGKPS09,DOURADO10a}; in view of these negative results, analyzing the behavior of these parameters in classes of graphs with an interesting structure, such as the class of Kneser graphs, is a natural research direction. Another objective of this work is to contribute to the knowledge on Kneser graphs. For instance, we characterize endpoints of diametral paths in Kneser graphs, and use this characterization as a step towards the determination of the geodetic number.

The remainder of this section provides all the necessary background. Section~\ref{sec:diam} gives a characterization of vertices that are endpoints of a diametral path in $K(2n+k,n)$, in terms of their intersection size. This characterization is used as a tool for proving Theorem~\ref{geoSet} on the geodetic number of $K(2n+k,n)$. Sections~\ref{sec:geodetic} and~\ref{sec:hull} present the main results on the geodetic number and the geodetic hull number of $K(2n+k,n)$, respectively. 

\bigskip

\noindent {\em Definitions and notation.}

\bigskip

Let $G =(V,E)$ be a finite graph. A {\em path} between $u,v\in V$ is a sequence of distinct vertices $u=x_1,x_2,\ldots,x_p=v$ such that $x_ix_{i+1}\in E$ for $1\leq i\leq p-1$, and its {\em length} is the number of edges therein. A {\em shortest path} between $u,v\in V$ is a path with minimum length. The {\em distance} between two vertices $u,v\in V$, denoted by $\mathit{dist}(u,v)$, is the length of a shortest path between them. The {\em diameter} of $G$ is defined as $\mathit{diam}(G)=\max_{u,v\in V} \mathit{dist}(u,v)$. A {\em diametral path} is a shortest path whose length is $\mathit{diam}(G)$. We denote the open neighborhood of a vertex $x$ by $N(x)$. For $S\subseteq V$, we define $N(S)=\cup_{x\in S} N(x)$.

An {\em $uv$-geodesic} is a shortest path between $u$ and $v$. The {\em geodetic interval} $I[u,v]$ is the set of all vertices belonging to some $uv$-geodesic. For a set $W\subseteq V(G)$, the geodetic interval $I[W]$ is defined as $I[W]=\cup_{u,v\in W} I[u,v]$. The set $W$ is {\em geodetically convex} (or {\em $g$-convex}) if $I[W]=W$, and a {\em geodetic set} of $G$ if $I[W]=V(G)$. The {\em geodetic number} $\mathit{gn}(G)$ of $G$ is the size of minimum geodetic set of $G$.

The family of $g$-convex sets of graph $G$ define the {\em geodetic convexity} associated with $G$. In general, a {\em convexity} associated with a graph $G$ consists of a collection $\mathcal{C}$ of subsets of $V(G)$, called {\em convex sets}, such that: (a) $\emptyset, V(G)\in\mathcal{C}$; (b) $\mathcal{C}$ is closed under intersections. Just like the geodetic convexity is defined over shortest paths, other graph convexities can be defined using different path systems, such as the {\em monophonic convexity}~\cite{DOURADO10}, associated with induced paths.

The {\em geodetic hull} of a set $W\subseteq V(G)$, denoted by $H[W]$, is the minimum $g$-convex set containing $W$. Also, $W$ is a {\em geodetic hull set} of $G$ if $H[W]=V(G)$. The {\em geodetic hull number}, denoted by $\mathit{ghn}\left(G\right)$, is the size of a minimum geodetic hull set of $G$. For an integer $k \geq 0$, we define $I^k[W]$ as follows: $I^0[W]=W$ and $I^{k}[W]=I[I^{k-1}[W]]$. It is not difficult to see that $H[W]=I^k[W]$ for some $k\geq 0$; in fact, $k$ can be taken as the minimum index for which $H[W]=I^k[W]$.

Let $n$ and $k$ be positive integers, and let $[n]=\{1,\ldots,n\}$ and $[2n+k]^n=\{S\subseteq [2n+k] : |S|=n\}$. The {\em Kneser graph} $K(2n+k,n)$ is the graph $G=(V,E)$ such that $V=[2n+k]^n$ and there is an edge between two vertices $u,v\in [2n+k]^n$ whenever $u\cap v=\emptyset$ (see~\cite{LOVASZ78}). Therefore, the Kneser graph $K(2n+k,n)$ contains $\binom{2n+k}{n}$ vertices and is a $\binom{n+k}{n}$-regular graph. The Kneser graph $K(5,2)$ is the well-known Petersen graph.

In~\cite{VALENCIAPABON2005} the authors show that $\mathit{diam}(K(2n+k,n))=\ceil{(n-1)/k}+1$. Additionally, for any $u,v\in V(K(2n+k,n))$ with $|u \cap v|=s$, they show that:

\begin{equation}\label{distInKneser}
\mathit{dist}(u,v) =
\left\{
	\begin{array}{cl}
		\min \{2\ceil{(n - s)/k}, 2\ceil{s/k} + 1\},  & \mbox{if } 1\leq k < n-1;\\
		2, & \mbox{if } k\geq  n-1.
	\end{array}
\right.
\end{equation}

Observe that if $n\geq 2$ and $k\geq n-1$ then ${\mathit diam}(K(2n+k,n))=2$. Thus, the graphs of diameter two form a infinite subfamily of ${\mathscr K}_n=\{K(2n+k,n) : k\geq 1\}$.

\section{Endpoints of diametral paths} \label{sec:diam}

The theorem below gives a necessary and sufficient condition for two vertices in the Kneser graph $K(2n+k,n)$ to be endpoints of a diametral path.

\begin{theorem}\label{intMax}
Let $u$ and $v$ be distinct vertices of $K(2n+k,n)$, and let $s=|u\cap v|$. Then $\mathit{dist}(u,v) = \mathit{diam}(K(2n+k,n))$ if and only if

\begin{equation}\label{boundsS}
\left(\left\lceil \frac{n-1}{2k} \right\rceil  - 1\right) k +1 \leq s \leq \left(\left\lceil \frac{n-1}{2k} \right\rceil - 1\right) k +1 + H(n,k),
\end{equation}
where
\begin{equation}\label{FunH}
H(n,k)=
\left\{
	\begin{array}{ll}
		 \max\{n\bmod k+k-2,0\},  & \mathrm{if} \ 0\leq n \bmod k\leq 1; \\
		 n\bmod k-2,  & \mathrm{if} \  2\leq n \bmod k\leq k-1. \\
	\end{array}
\right.
\end{equation}

\end{theorem}

\noindent{\bf Proof:} First, assume condition (\ref{boundsS}) holds.

If $k\geq n-1$, we have $\mathit{diam}(K(2n+k,n))=2$ and $s\geq 1$, that is, $u$ and $v$ are not adjacent. Therefore, $\mathit{dist}(u,v)=2=\mathit{diam}(K(2n+k,n))$.

Now, consider $n-1 = 2\alpha k+\beta k + \gamma$ for $\alpha \in \mathbb{N}_0$, $\beta \in \{0,1\}$, and $0\leq \gamma<k$. For $s=\left(\left\lceil \frac{n-1}{2k} \right\rceil - 1\right) k +1+\Delta$ with $0\leq \Delta\leq H(n,k)$, we analyze two cases:

\medskip

\noindent {\em Case 1:} \ $\beta=\gamma=0$.

\medskip

In this case, $n=2\alpha k+1$, $s=(\alpha-1)k+1+\Delta$, and $\mathit{diam}(K(2n+k,n))=\ceil{(n-1)/k}+1=2\alpha +1$.

If $k=1$, we have $H(n,k)=\Delta=0$. By using Eq.~(\ref{distInKneser}), it follows that
\begin{flalign*}
\mathit{dist}(u,v) &= \min \{2\ceil{((\alpha+1)k-\Delta)/k}, 2\ceil{((\alpha -1)k+1+\Delta)/k} + 1\}\\
&= \min \{2\ceil{(\alpha+1)/1}, 2\ceil{(\alpha -1)+1/1} + 1\}\\
&= 2\alpha +1 = \mathit{diam}(K(2n+k,n)).
\end{flalign*}

Let us assume now $k\geq 2$. Since $n\bmod k=1$, we have $H(n,k)=k-1$, and, thus, $0\leq \Delta\leq k-1$. Then, by Eq.~(\ref{distInKneser}), it follows that
\begin{flalign*}
\mathit{dist}(u,v) &= \min \{2\ceil{((\alpha+1)k-\Delta)/k}, 2\ceil{((\alpha -1)k+1+\Delta)/k} + 1\}\\
&= \min \{2\ceil{(\alpha+1)-\Delta/k}, 2\ceil{(\alpha -1)+(1+\Delta)/k} + 1\}\\
&= \min \{2\alpha+2, 2\alpha + 1\}\\
&= 2\alpha + 1 = \mathit{diam}(K(2n+k,n)),
\end{flalign*}

\noindent for every $0\leq \Delta\leq H(n,k)=k-1$. This concludes Case 1.

\medskip

\noindent {\em Case 2:} \ $\beta= 1$ or $\gamma\neq 0$.

\medskip

In this case, we have $s=\alpha k+1+\Delta$.

First, let us assume $k=1$. Then, $\gamma=\Delta=H(n,k)=0$, which implies $\beta=1$ and $\mathit{diam}(K(2n+k,n))=\ceil{(2\alpha+1)/1}+1=2\alpha+2$. By using Eq.~(\ref{distInKneser}) and the fact that $n=2\alpha k+\beta k+\gamma+1$, we have that:
\begin{flalign}
\mathit{dist}(u,v) &= \min \{2\ceil{(\alpha k+ \beta k + \gamma -\Delta)/k}, 2\ceil{(\alpha k+ 1 +\Delta)/k} + 1\}\nonumber\\
&= \min \{2\ceil{\alpha +\beta  + (\gamma -\Delta)/k}, 2\ceil{\alpha +(1+\Delta)/k} + 1\}\label{distCase2}\\
&= \min \{2\ceil{\alpha+1}, 2\ceil{\alpha +1} + 1\}\nonumber\\
&= 2\alpha+2 = \mathit{diam}(K(2n+k,n)).\nonumber
\end{flalign}

Assume now $k\geq 2$. Note that $n\bmod k =(\gamma+1)\bmod k$. Sub\-sti\-tut\-ing in Eq.~(\ref{FunH}):
\begin{equation*}
    H(n,k) =
\left\{
	\begin{array}{ll}
		\max\{(\gamma+1)\bmod k+k-2,0\}, &\mbox{if } 0\leq(\gamma+1)\bmod k\leq 1;\\
		(\gamma+1)\bmod k-2,  &\mbox{if }  2\leq(\gamma+1)\bmod k\leq k-1.\\
	\end{array}
\right.
\end{equation*}

From the above equation, $\gamma-H(n,k)=1-k$ if $\gamma=0$, and $\gamma-H(n,k)=1$ if $1\leq\gamma\leq k-1$. Then, we have  $1\leq\gamma-\Delta\leq k-1$ when $\gamma\neq 0$. We analyze the possible cases in Eq.~(\ref{distInKneser}):

\begin{itemize}
\item[--] If $\beta = 1$ and $\gamma\neq 0$, $\mathit{dist}(u,v)=\min\{2\alpha+4,2\alpha+3\}=2\alpha + 3$ and $\mathit{diam}(K(2n+k,n))=\ceil{(2\alpha k+\beta k +\gamma)/k}+1=2\alpha+3$;
\item[--] If $\beta=0$ and $\gamma\neq0$, $\mathit{dist}(u,v)=\min\{2\alpha+2, 2\alpha+3\}=2\alpha+2=\mathit{diam}(K(2n+k,n))$;
\item[--] If $\beta=1$ and $\gamma=0$, $\mathit{dist}(u,v)=\min\{2\alpha+2,2\alpha+3\}=2\alpha+2$ and $\mathit{diam}(K(2n+k,n))=\ceil{(2\alpha k+\beta k)/k}+1=2\alpha+2$.
\end{itemize}

This concludes Case 2 and the first part of the proof.

Conversely, suppose $\mathit{dist}(u,v)=\mathit{diam}(K(2n+k,n))$, and assume by contradiction that there is an integer $\varepsilon\geq 1$ such that
$$\left(\left\lceil\frac{n-1}{2k}\right\rceil-1\right)k+1-\varepsilon=s< \left(\left\lceil\frac{n-1}{2k}\right\rceil-1\right)k+1 \ \ \      \text{or}
$$
$$
\left(\left\lceil\frac{n-1}{2k}\right\rceil-1\right)k+1+H(n,k)>s=\left(\left\lceil\frac{n-1}{2k}\right\rceil-1\right)k+1+H(n,k)+\varepsilon.$$

Assuming $n-1 = 2\alpha k+\beta k + \gamma$ for $\alpha \in \mathbb{N}_0$, $\beta \in \{0,1\}$ again, we analyze the two possible cases.

\medskip

\noindent {\em Case 1:} \
$s=\left(\left\lceil\frac{n-1}{2k}\right\rceil-1\right) k+1-\varepsilon$.

\medskip

We divide the proof of Case 1 in two subcases, analyzing possible values of $\beta$ and $\gamma$.

\medskip

\noindent {\em Case 1.1:} \ $\beta=\gamma=0$.

\medskip

This case implies $n=2\alpha k+1$ and $s=(\alpha -1)k+1-\varepsilon$. Thus, substituting in Eq.~(\ref{distInKneser}),
\begin{flalign*}
\mathit{dist}(u,v) &= \min \{2\ceil{((\alpha+1)k+\varepsilon)/k}, 2\ceil{((\alpha -1)k+1-\varepsilon)/k} + 1\}\\
&= \min \{2\ceil{(\alpha+1)+\varepsilon/k}, 2\ceil{(\alpha -1)+(1-\varepsilon)/k} + 1\}\\
&= 2\ceil{(\alpha -1)+(1-\varepsilon)/k} + 1\leq 2\alpha -1.
\end{flalign*}
However, $\mathit{diam}(K(2n+k,n))=\ceil{(n-1)/k}+1=\ceil{2\alpha k/k}+1=2\alpha+1$.

\medskip

\noindent {\em Case 1.2:} \ $\beta=1$ or $\gamma\neq 0$.

\medskip

This case implies $n=2\alpha k+\beta k+\gamma+1$ and $s=\alpha k+1-\varepsilon$. Substituting in Eq.~(\ref{distInKneser}),
\begin{flalign*}
\mathit{dist}(u,v) &= \min\{2\ceil{(\alpha k+\beta k +\gamma+\varepsilon)/k}, 2\ceil{(\alpha k+1-\varepsilon)/k}+1\}\\
&= \min \{2\ceil{\alpha+\beta+(\gamma+\varepsilon)/k}, 2\ceil{\alpha +(1-\varepsilon)/k}+1\} \\
&= 2\ceil{\alpha+(1-\varepsilon)/k}+1\leq 2\alpha +1.
\end{flalign*}
However, for $\beta=1$ or $\gamma\neq 0$, $2\alpha+2\leq\mathit{diam}(K(2n+k,n))\leq2\alpha+3$. Therefore, both Cases 1.1 and 1.2 lead to contradictions, and this concludes Case 1.

\medskip

\noindent {\em Case 2:} \
$s=\left(\left\lceil\frac{n-1}{2k}\right\rceil-1\right) k+1+H(n,k)+\varepsilon$.

\medskip

Again, we analyze the possible values of $\beta$ and $\gamma$.

\medskip

\noindent {\em Case 2.1:} \ $\beta=\gamma=0$.

\medskip

This case implies $n=2\alpha k+1$ and $s=(\alpha-1)k+1+H(n,k)+\varepsilon$.

Assume $k=1$. Then, $H(n,k)=0$, and using Eq.~(\ref{distInKneser}),
\begin{flalign*}
\mathit{dist}(u,v) &= \min\{2\ceil{((\alpha+1)k-\varepsilon)/k}, 2\ceil{((\alpha -1)k+1+\varepsilon)/k}+1\}\\
&= \min\{2(\alpha+1-\varepsilon),2\alpha+2\varepsilon+1\}\\
&= 2(\alpha+1-\varepsilon)\leq 2\alpha,
\end{flalign*}
\noindent while $\mathit{diam}(K(2n+k,n))=\ceil{(n-1)/k}+1=\ceil{2\alpha k/k}+1=2\alpha+1$.

Assume now $k\geq 2$. Then, recall that $H(n,k)=k-1$, since $n\bmod k=1$. Substituting in Eq.~(\ref{distInKneser}), we have that
\begin{flalign*}
\mathit{dist}(u,v) &= \min\{2\ceil{((\alpha+1)k-H(n,k)-\varepsilon)/k}, 2\ceil{((\alpha-1)k+1+H(n,k)+\varepsilon)/k}+1\}\\
&= \min\{2\ceil{(\alpha+1)-(H(n,k)+\varepsilon)/k}, 2\ceil{(\alpha-1)+(1+H(n,k)+\varepsilon)/k}+1\}\\
&= \min\{2\ceil{(\alpha+1)-(k-1+\varepsilon)/k}, 2\ceil{(\alpha -1)+(k+\varepsilon)/k}+1\}\\
&= 2\ceil{(\alpha+1)-(k-1+\varepsilon)/k}
\leq2\alpha<\mathit{diam}(K(2n+k,n))=2\alpha+1.
\end{flalign*}

\medskip

\noindent {\em Case 2.2:} \ $\beta=1$ or $\gamma\neq 0$.

\medskip

In this case, $n=2\alpha k+\beta k+\gamma+1$ and $s=\alpha k+1+H(n,k)+\varepsilon$.

If $k=1$, we have $H(n,k)=\gamma=0$, $\beta=1$, and $\mathit{diam}(K(2n+k,n))=2\alpha+2$.
Thus, by using Eq.~(\ref{distInKneser}), we have that
\begin{flalign*}
\mathit{dist}(u,v) &= \min\{2\ceil{(\alpha k+\beta k+\gamma-H(n,k)-\varepsilon)/k}, 2\ceil{(\alpha k+1+H(n,k)+\varepsilon)/k}+1\}\\
&= \min\{2\ceil{\alpha+1-\varepsilon/k}, 2\ceil{\alpha+(1+\varepsilon)/k}+1\}\\
&= 2\ceil{\alpha+1-\varepsilon/k}
\leq2\alpha<\mathit{diam}(K(2n+k,n))=2\alpha+2.
\end{flalign*}

Finally, assume $k\geq 2$. Again, by Eq.~(\ref{distInKneser}):
\begin{flalign}
\mathit{dist}(u,v) &= \min\{2\ceil{(\alpha k+\beta k+\gamma-H(n,k)-\varepsilon)/k},
2\ceil{(\alpha k+1+H(n,k)+\varepsilon)/k}+1\}\nonumber\\
&= \min \{2\ceil{\alpha+\beta+(\gamma-H(n,k)-\varepsilon)/k}, 2\ceil{\alpha +(1+H(n,k)+\varepsilon)/k}+1\}\label{eq:lastcase}.
\end{flalign}

Recall from Case 2 in the first part of the proof that $\gamma-H(n,k)=1$, for all $1\leq\gamma\leq k-1$. By confronting Eq.~(\ref{eq:lastcase}) with the possibilities for $\beta$ and $\gamma$, we have:

\begin{itemize}
\item[--] If $\beta = 1$ and $\gamma\neq 0$,
    \begin{flalign*}
    \mathit{dist}(u,v) &= \min\{2\ceil{\alpha+1+(1-\varepsilon)/k}, 2\ceil{\alpha+(1+H(n,k)+\varepsilon)/k}+1\}\\
    &= 2\ceil{\alpha+1+(1-\varepsilon)/k}
    \leq 2\alpha+2<\mathit{diam}(K(2n+k,n))=2\alpha+3.
    \end{flalign*}

\item[--] if $\beta=0$ and $\gamma\neq0$,
    \begin{flalign*}
    \mathit{dist}(u,v) &= \min\{2\ceil{\alpha+(1-\varepsilon)/k}, 2\ceil{\alpha+(1+H(n,k)+\varepsilon)/k}+1\}\\
    &= 2\ceil{\alpha+(1-\varepsilon)/k}\leq 2\alpha<\mathit{diam}(K(2n+k,n))=2\alpha+2.
    \end{flalign*}

\item[--] If $\beta=1$ and $\gamma=0$, we have $-(H(n,k)+\varepsilon)\leq -1$. Thus,
    \begin{flalign*}
    \mathit{dist}(u,v) &=\min\{2\ceil{\alpha +1-(H(n,k)+\varepsilon)/k}, 2\ceil{\alpha +(1+H(n,k)+\varepsilon)/k}+1\}\\
    &= 2\ceil{\alpha+1-(H(n,k)+\varepsilon)/k}
    \leq 2\alpha<\mathit{diam}(K(2n+k,n))=2\alpha+2.
    \end{flalign*}
\end{itemize}

Therefore, both Cases 2.1 and 2.2 lead to contradictions. This concludes Case 2 and the proof of the theorem. \ $\Box$

\medskip

An interesting fact derived from Theorem~\ref{intMax} is:

\begin{corollary}\label{intProp}
Let $u,v\in V(K(2n+k,n))$. Then:
\begin{itemize}
    \item if $|u\cap v|<\left(\left\lceil\frac{n-1}{2k}\right\rceil-1 \right) k+1$ then $\mathit{dist}(u,v)$ is odd;
    \item if $|u\cap v|>\left(\left\lceil\frac{n-1}{2k}\right\rceil-1 \right) k+1+H(n,k)$ then $\mathit{dist}(u,v)$ is even.
\end{itemize}
\end{corollary}

\noindent{\bf Proof:}
According to the proof of Theorem~\ref{intMax}, Cases 1 and 2 guarantee that intersections with fewer than $\left(\left\lceil\frac{n-1}{2k}\right\rceil-1\right) k+1$ elements imply odd distances, while intersections with more than $\left(\left\lceil\frac{n-1}{2k}\right\rceil-1\right) k+1+H(n,k)$ imply even distances. \ $\Box$

\begin{figure}[htbp]
\begin{center}
\includegraphics[scale=0.35]{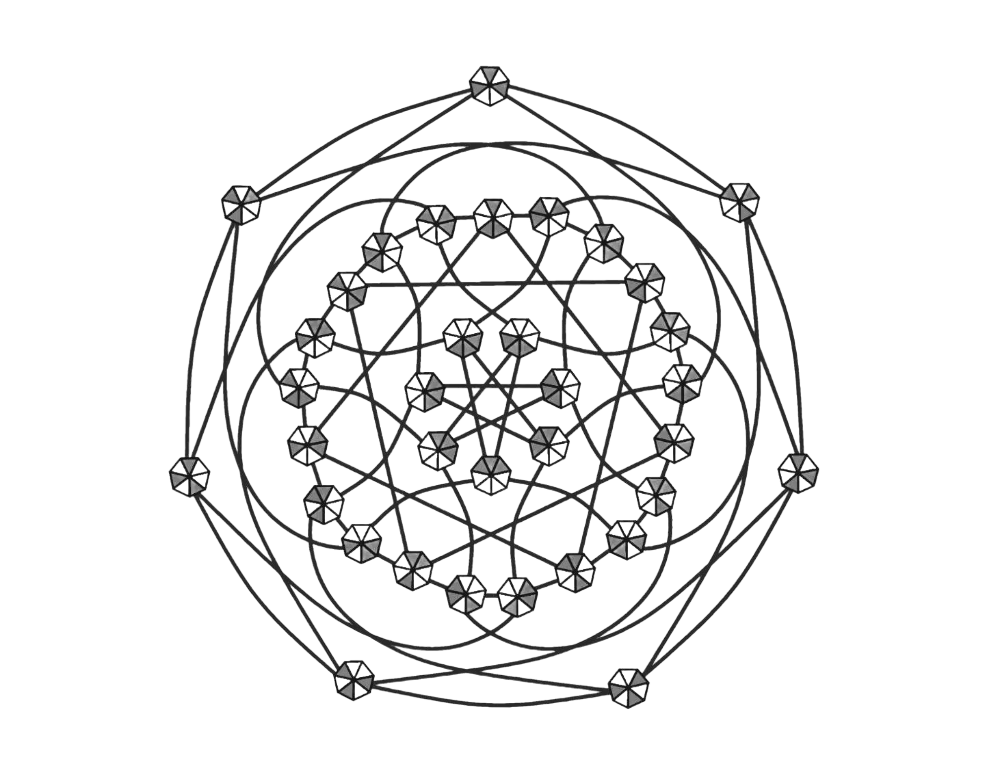}
\caption{Kneser graph $K(7,3)$ ($n=3$ and $k=1$).}\label{fig:kneser73}
\end{center}
\end{figure}

\begin{example}{\em
Figure~\ref{fig:kneser73} depicts the Kneser graph $K(7,3)$, for which $n=3$ and $k=1$. By inspection (or using the formula in~\cite{VALENCIAPABON2005}), $\mathit{diam}(K(7,3))=3$. Substituting the values of $n$ and $k$ in Eqs.~(\ref{boundsS}) and~(\ref{FunH}), two vertices $u$ and $v$ are endpoints of a diametral path in $K(7,3)$ if and only if $|u\cap v|=1$. For instance, vertices $\{1,2,3\}$ and $\{1,4,5\}$ are endpoints of the diametral path $P=\{1,2,3\},\{4,5,6\},\{2,3,7\},\{1,4,5\}$. By Corollary~\ref{intProp}, if $|u\cap v|=0<1$ then $u$ and $v$ are at an odd distance $d=1$, and if $|u\cap v|=2>1$ then $u$ and $v$ are at an even distance $d=2$. For instance, vertices $\{1,2,3\}$ and $\{1,2,4\}$ are endpoints of the path $P'=\{1,2,3\},\{5,6,7\},\{1,2,4\}$.          }
\end{example}

\section{Geodetic number}\label{sec:geodetic}

The next theorem gives a sufficient condition for a set $S\subseteq V(K(2n+k,n))$ to be a geodetic set. Say that two vertices $u,v$ in a graph $G$ are {\em diametrically opposed} if $d(u,v)={\mathit diam}(G)$.

\begin{theorem}\label{geoSet}
Let $r\in V(K(2n+k,n))$ and let $D$ be the set of all vertices of  $K(2n+k,n)$ diametrically opposed to $r$. Then $D\cup\{r\}$ is a geodetic set.
\end{theorem}

\noindent\textbf{Proof:} \ Let $T(r)$ be a tree rooted at $r$, obtained by a breadth-first search in $K(2n+k,n)$. Let $L(x)$ be the level of a vertex $x$ in $T(r)$. Trivially, $L(r)=0$. In addition, $x\in D$ if and only if $L(x)=\ceil{(n-1)/k}+1$. Let $\Le(i)=\{x : L(x)=i\}$.

In order to prove that $D\cup\{r\}$ is a geodetic set, we show that each $x\in\Le(i)$, with $0\leq i<\ceil{(n-1)/k}+1$, has at least one neighbor $z\in\Le(i+1)$. This is trivial for $i=0$. By Theorem~\ref{intMax}, either $|r\cap x|<\left(\left\lceil \frac{n-1}{2k} \right\rceil-1\right) k+1$ or $|r\cap x|>\left(\left\lceil\frac{n-1}{2k} \right\rceil-1\right) k+1+H(n,k)$.

Let $y\in N(x)\cap\Le(i-1)$. We analyze two cases.

\medskip

\noindent {\em Case 1: $i$ is even.}

\medskip

In this case, by Corollary~\ref{intProp}, $|r\cap y|<\left(\left\lceil\frac{n-1}{2k}\right\rceil-1\right) k+1$.

Let $k'=\min\{k,\left(\left\lceil\frac{n-1}{2k}\right\rceil-1\right) k+1-|r\cap y|\}$, and let $z_1$ be a set formed by $k'$ elements in $r\setminus(x\cup y)$. Observe that $z_1$ exists because ${\mathit dist}(x,r) <{\mathit diam}(K(2n+k,n))$, which in turn implies the existence of a set that has at least $k'$ more  elements in common with $r$ than $y$. Additionally, let $z_2$ be a set formed by $n-k'-|r\cap y|$ elements in $y\setminus r$. It follows that $z=(r\cap y)\cup z_1\cup z_2$ is a neighbor of $x$ such that $L(z)=i+1$, and this concludes Case 1.

\medskip

\noindent {\em Case 2: $i$ is odd.}

\medskip

In this case, by Corollary~\ref{intProp}, $|r\cap y|>\left(\left\lceil \frac{n-1}{2k}\right\rceil-1\right) k+1+H(n,k)$.

Let $k''=\min\{k,|r\cap y|-(\left(\left\lceil \frac{n-1}{2k}\right\rceil -1\right) k+1+H(n,k))\}$, and let $z_3$ be a set formed by $k''$ elements in $r\cap y$. In addition, let $z_4$ be formed by $n-k''$ elements in $[2n+k]\setminus(r\cup x)$. Again, $z_3$ and $z_4$ exist because there is a set with at least $k''$ elements that intersect $r$ less than $y$, since ${\mathit dist}(x,r) <{\mathit diam}(K(2n+k,n))$. To conclude Case 2, note that $z=(y\setminus z_3) \cup z_4$ is a neighbor of $x$ with $L(z)=i+1$.

\medskip

As seen above, each vertex in $\Le(\ceil{(n-1)/k})$ has a neighbor in $D=\Le(\ceil{(n-1)/k}+1)$. Therefore, $x\in I[r,x']$, for some $x'\in D$. In other words, $D\cup\{r\}$ is a geodetic set. \ $\Box$

\begin{corollary}\label{coro:bound-gn}
Let $p=\left(\left\lceil\frac{n-1}{2k}\right\rceil-1\right) k+1$. Then,
\begin{equation}\label{eq:bound-gn}
    \mathit{gn}(K(2n+k,n))\leq 1+\sum_{i=p}^{p+H(n,k)}\binom{n}{i} \binom{n+k}{n-i}.
\end{equation}
\end{corollary}

\noindent\textbf{Proof:} \ The bound in Eq.~\ref{eq:bound-gn} is precisely the size of $D\cup\{r\}$ in the proof of Theorem~\ref{geoSet}. \ $\Box$

\medskip

Note that the bound in Eq.~\ref{eq:bound-gn} is valid for all possible diameter values of $K(2n+k,n)$. If $\mathit{diam}(K(2n+k,n))=2$, we can improve the result of Theorem~\ref{geoSet} and find the exact value of $\mathit{gn}(K(2n+k,n))$, as shown in the next theorem:

\begin{theorem}\label{theo:3.2}
If $k\geq n-1$ then
\begin{equation*}
\mathit{gn}(K(2n+k,n))= \binom{2n+k-1}{n-1}.
\end{equation*}
\end{theorem}

\noindent\textbf{Proof:} \ Let $u,v\in V(K(2n+k,n))$ and $s=|u\cap v|$. Since ${\mathit diam}(K(2n+k,n))=2$, the elements of a geodetic set of $K(2n+k,n)$ must have $1\leq s\leq n-1$.

We notice that, for any geodetic set $X$ and a vertex $x\in V(K(2n+k,n))$, if $x$ has no neighbors in $X$ then $x\in X$, because the diameter is two. We construct a set $D\subseteq V(K(2n+k,n))$ consisting of pairwise diametrically opposed vertices as follows: for each $s=1,\ldots,n-1$, include in $D$ a maximal subset $D_s$ such that, for every distinct $u,v\in D_s$, $u\,\cap\,v = \{1,2,\ldots,s\}$. In other words, $D_1$ is formed by vertices with pairwise intersection $\{1\}$, $D_2$ by vertices with pairwise intersection $\{1,2\}$, and so on. Note that $D$ is an independent set. Moreover, one can verify that
$$|D|=  \binom{2n+k-1}{n-1}.$$


The above construction generates a set $D$ of maximal size in which any pair of vertices in $D$ has a diametral path connecting them. Moreover, notice that any subset $S$ with $|D|-1$ vertices implies some vertex with no edge to vertices in $S$. Then, $|D|\leq |S|$ for any geodetic set $S$. Now, let $x\in V(K(2n+k,n))\setminus D$. Then $1\notin x$, and, by construction, $x$ has distinct neighbors $y,z\in D$ (recall that $y\cap z\neq\emptyset$, which implies $|x\cup y\cup z|\leq 3n-1\leq 2n+k$). Thus, $D$ is a geodetic set of minimum size and $\mathit{gn}(K(2n+k,n))=|D|$. \ $\Box$

\begin{example}\label{ex:kneser62}
{\em Let $n=k=2$, and consider the Kneser graph $K(6,2)$. According to the proof of Theorem~\ref{geoSet}, let $r=\{1,2\}$ and let $D$ be the following set:
$$D=\{\{x,3\},\{x,4\},\{x,5\},\{x,6\} : x\in\{1,2\}\}.$$
Note that $D$ is the set of vertices diametrically opposed to $r$. Thus, $D\cup\{r\}$ is a geodetic set of $K(6,2)$, with size $9$. Indeed, substituting $n=k=2$ in Eq.~\ref{eq:bound-gn} we have $p=1$, $H(n,k)=0$, and
$$\mathit{gn}(K(6,2))\leq 1+\binom{2}{1} \binom{4}{1}=9.$$
However, we can improve this result using Theorem~\ref{theo:3.2}, whose proof tells us that a maximal set $D'$ of vertices with pairwise intersection $I=\{1\}$ is a minimum geodetic set of $K(6,2)$. Thus, $D'=\{ \{1,2\}, \{1,3\}, \{1,4\}, \{1,5\}, \{1,6\}\}$ is the required set and $\mathit{gn}(K(6,2))=5$. Note that $D'$ is also an independent set.}
\end{example}


\section{Geodetic hull number}\label{sec:hull}

The next lemma gives a necessary condition for a set to be a geodetic hull set of a Kneser graph with diameter two.

\begin{lemma}\label{lem:hull}
Suppose $k\geq n-1$, and let $x,y,z\in V(K(2n+k,n))$ be vertices such that 
$$(x\cup y\cup z)\setminus(x\cap y\cap z)=\{h,i,j\},$$
where $h\in x$, $i\in y$, and $j\in z$. Then, $\{x,y,z\}$ is a geodetic hull set of $K(2n+k,n)$.
\end{lemma}

\noindent \textbf{Proof:} \ Let $w\in V(K(2n+k,n))\setminus\{x,y,z\}$. If $w$ has empty intersections with at least two of $x,y,z$, say $x$ and $y$, then $w\in I[x,y]$. Then we may assume, without loss of generality, that $w\cap x\neq\emptyset$ and $w\cap y\neq\emptyset$.

If $|w\cap(x\cup y)|>1$ and $w\cap z=\emptyset$ then $w\cap x\cap y =\{h,i\}$. In addition, we can define a vertex $w'$ consisting of element $j$ plus $n-1$ elements of $[2n+k]\setminus(w\cup x\cup y)$. Notice that $w'\in I[x,y]$ and, consequently, $w\in I[w',z]$.

If $|w\cap(x\cup y)|=1$ and $w\cap z= \emptyset$ then either $w\cap(x\cup y)=\{x\}$ or $w\cap(x\cup y)=\{y\}$. The former case implies $w\in I[y,z]$, and the latter $w\in I[x,z]$.
%
%
Such observations imply that $N(z)\subseteq H[\{x,y,z\}]$. Likewise, if $w$ has nonempty intersections with $x\cup z$ (resp., $y\cup z$) then  $N(y)\subseteq H[\{x,y,z\}]$ (resp., $N(x)\subseteq H[\{x,y,z\}]$). Therefore, $N(\{x,y,z\})\subseteq H[\{x,y,z\}]$.

Finally, if $w\notin N(\{x,y,z\})$, then $\mathit{dist}(w,a)=2$ for $a\in\{x,y,z\}$. This means that $w\in I[b,c]$ for distinct $b,c\in N(\{x,y,z\})$. Hence, $\{x,y,z\}$ is a geodetic hull set. \ $\Box$

\begin{example}
{\em Consider again the Kneser graph $K(6,2)$, and let $x=\{1,2\}$, $y=\{1,3\}$, and $z=\{1,4\}$. Note that $(x\cup y\cup z)\setminus(x\cap y\cap z)=\{2,3,4\}$, with $2\in x$, $3\in y$, and $4\in z$. Therefore, according to Lemma~\ref{lem:hull}, $S=\{x,y,z\}$ is a geodetic hull set of $K(6,2)$. Indeed,
$$I^1[\{x,y,z\}] = \{x,y,z\}\cup \{\{4,5\},\{4,6\},\{5,6\},\{3,5\},\{3,6\},\{2,5\},\{2,6\}\}$$
$$\text{and}$$
$$I^2[\{x,y,z\}] = I^1[\{x,y,z\}]\cup \{\{1,5\},\{1,6\},\{2,3\},\{2,4\},\{3,4\}\} = V(K(6,2)),$$
that is, $H[\{x,y,z\}] = V(K(6,2))$.}
\end{example}

\begin{theorem}\label{theo:4.1}
If $k\geq n-1$, then
\begin{equation}
\mathit{ghn}(K(2n+k,n)) =
\left\{
	\begin{array}{cl}
		2, & \mbox{if } k>2;\\
		3, & \mbox{otherwise.}
	\end{array}
\right.
\end{equation}
\end{theorem}

\noindent \textbf{Proof:} \  Recall that if $k\geq n-1$, then ${\mathit diam}(K(2n+k,n))=2$. Suppose $k>2$, and let $x,y,z\in V(K(2n+k,n))$ such that $(x\cup y\cup z)\setminus(x\cap y\cap z)=\{h,i,j\}$, where $h\in x$, $i\in y$, and $j\in z$. By Lemma~\ref{lem:hull}, we know that $\{x,y,z\}$ is a geodetic hull set. We show that $\{x,y\}$ is still a geodetic hull set in this case.

Notice that there are $w,w'\in I[x,y]$ such that $|w\cup w'|=n+1$ and $j\notin w\cup w'$. Also, $|x\cup y\cup z\cup w\cup w'|=|x\cup y\cup z|+|w\cup w'|=(n+2)+(n+1)=2n+3$ and $z\in I[w,w']$. Therefore, $z\in I^2[\{x,y\}]$ and this concludes the case $k>2$.

Now, suppose $k\leq 2$. In this case, we show that no set $S$ with $|S|=2$ is a geodetic hull set. Let $S=\{x',y'\}$, and assume $|x'\cap y'|=s$. Notice that $I[x',y']$ must contain at least two vertices $w$ and $w'$, since they must provide vertices with no edge to $x'$ or $y'$. Observe that, even for $|w\cap w'|=1$,  $|x'\cup y'\cup w\cup w'|=|x'\cup y'|+|w\cup w'|=(2n-s+1)+(n+1)=3n-s+2$. But since $2n+k\leq 2n+2$, we have $3n-s+2\leq 2n+2$. This implies $s\geq n$, a contradiction. In other words, $H[S]\neq V(K(2n+k,n))$. Thus, if $k\leq 2$, $\mathit{ghn}(K(2n+k,n))=3$. This concludes the proof of the theorem. \ $\Box$

\begin{corollary}\label{coro:diam2-hn3}
The only Kneser graphs with diameter two and geodetic hull number three are $K(5,2)$, $K(6,2)$, and $K(8,2)$.
\end{corollary}

\noindent \textbf{Proof:} \ By Theorem~\ref{theo:4.1}, the Kneser graphs with diameter two and geodetic hull number three are obtained by combining the inequalities $k\geq n-1\geq 1$ and $k\leq 2$. \ $\Box$





\end{document}